\documentclass[reqno]{article}%
\usepackage[intlimits]{amsmath}
\usepackage{graphicx}
\usepackage{amsfonts}
\usepackage{amssymb}%
\newtheorem{theorem}{Theorem}

\newtheorem{corollary}[theorem]{Corollary}

\newtheorem{lemma}[theorem]{Lemma}

\newenvironment{proof}[1][Proof]{\noindent{\textbf {#1}  }}  {\hfill$\Box$}

\begin{document}

\title{The Asymptotics of Strongly Regular Graphs}
\author{V. Nikiforov\\{\small Department of Mathematical Sciences}\\{\small University of Memphis, Memphis, TN 38152}}
\maketitle

\begin{abstract}
A strongly regular graph is called trivial if it or its complement is a union
of disjoint cliques. We prove that the parameters $n,k,\lambda,\mu$ of
nontrivial strongly regular graphs satisfy%
\[
\lambda=k^{2}/n+o\left(  n\right)  \text{ \ \ and \ \ }\mu=k^{2}/n+o\left(
n\right)  .
\]

It follows, in particular, that every infinite family of nontrivial strongly
regular graphs is quasi-random in the sense of Chung, Graham and Wilson.

\end{abstract}

\section{Introduction}

Our graph-theoretic notation is standard (see, e.g. \cite{Bol98}). Given a
graph $G$ and a set $R\subset V\left(  G\right)  ,$ we write $\widehat
{d}\left(  R\right)  $ for the number vertices in $G$ joined to every vertex
in $R$ and call the value $\widehat{d}\left(  R\right)  $ the \emph{codegree}
of $R$.

A \emph{strongly regular graph} (\emph{srg} for short) with parameters
$n,k,\lambda,\mu$ is a $k$-regular graph of order $n$ such that $\widehat
{d}\left(  uv\right)  =\lambda$ if $uv$ is an edge, and $\widehat{d}\left(
uv\right)  =\mu$ if $uv$ is not an edge; we denote by $SR\left(
n,k,\lambda,\mu\right)  $ a srg with parameters $n,k,\lambda,\mu$.

Observe that any graph $rK_{m}$ is an $SR\left(  mr,m-1,m-2,0\right)  ;$ we
call these graphs and their complements \emph{trivial }srgs.

Srgs have been intensively studied; we refer the reader to, e.g.
\cite{CavL91}, \cite{Bro96}, and \cite{Cam04}. Among the many problems related
to srgs, probably the most intriguing one is to find strong necessary
conditions for the parameters of a srg. Despite the numerous partial results,
no exact condition of wide scope is known. If we look for asymptotic
conditions, however, the problem becomes more tangible.

In this note we investigate the parameters of nontrivial srgs when the order
tends to infinity. Somewhat surprisingly it turns out that the parameters
$\lambda$ and $\mu$ are asymptotically equal. More precisely, the following
theorem holds.

\begin{theorem}
\label{main} The parameters $n,k,\lambda,\mu$ of nontrivial strongly regular
graphs satisfy
\begin{equation}
\lambda=k^{2}/n+o\left(  n\right)  \text{ \ \ and \ \ }\mu=k^{2}/n+o\left(
n\right)  . \label{cdc}%
\end{equation}

\end{theorem}

In terms of quasi-random graphs (e.g., see \cite{CGW89}, \cite{KrSu}) this
result implies that every infinite family of nontrivial srgs is quasi-random.

Recently Cameron \cite{Cam03} discussed the randomness aspect of srgs;
however, already Thomason \cite{Tho87} suggested that close relations between
srgs and quasi-random graphs might exist. Our result shows that, in fact,
there is a straightforward relationship.

To prove Theorem \ref{main} we shall use Semer\'{e}di's Uniformity Lemma (SUL
for short) - a widely applicable tool in extremal graph theory, but seldom, if
ever, applied to \textquotedblleft rigid\textquotedblright\ combinatorial
objects like srgs.

In Section \ref{sul} we give the notions related to SUL and several counting
lemmas; the proof of Theorem \ref{main} is presented in Section \ref{prf}.

\section{\label{sul}Semer\'{e}di's Uniformity Lemma}

For expository matter on Szemer\'{e}di's uniformity lemma (SUL) the reader is
referred to \cite{KoSi93} and \cite{Bol98}. This remarkable result is usually
called Szemer\'{e}di's Regularity Lemma, but the term \textquotedblleft
uniformity\textquotedblright\ seems more appropriate to its spirit.

We shall introduce some notation. Given a graph $G,$ if $u\in V\left(
G\right)  $ and $Y\subset V\left(  G\right)  ,$ we write $d_{Y}\left(
u\right)  $ for the number of neighbors of $u$ in $Y;$ similarly, if $R\subset
V\left(  G\right)  ,$ we write $\widehat{d}_{Y}\left(  R\right)  $ for the
number of vertices in $Y$ that are joined to every vertex in $R.$ The set of
neighbors of a vertex $u$ is denoted by $\Gamma\left(  u\right)  $.

Let $G$ be a graph; if $A,B\subset V\left(  G\right)  $ are nonempty disjoint
sets, we write $e\left(  A,B\right)  $ for the number of $A-B$ edges; the
value
\[
d\left(  A,B\right)  =\frac{e\left(  A,B\right)  }{\left\vert A\right\vert
\left\vert B\right\vert }%
\]
is called the \emph{density} of the pair $\left(  A,B\right)  .$

Let $\varepsilon>0;$ a pair $\left(  A,B\right)  $ of two nonempty disjoint
sets $A,B\subset V\left(  G\right)  $ is called $\varepsilon$\emph{-uniform}
if the inequality%
\[
\left\vert d\left(  A,B\right)  -d\left(  X,Y\right)  \right\vert <\varepsilon
\]
holds for every $X\subset A,$ $Y\subset B$ with $\left\vert X\right\vert
\geq\varepsilon\left\vert A\right\vert $ and $\left\vert Y\right\vert
\geq\varepsilon\left\vert B\right\vert .$

We shall use SUL in the following form.

\begin{theorem}
[Szemer\'{e}di's Uniformity Lemma]\label{SUL} Let $l\geq1$, $\varepsilon>0$.
There exists $M=M\left(  \varepsilon,l\right)  $ such that, for every graph
$G$ of sufficiently large order $n$, there exists a partition $V\left(
G\right)  =\cup_{i=0}^{p}V_{i}$ satisfying $l\leq p\leq M$ and:

$\emph{(i)}$ $\left\vert V_{0}\right\vert <\varepsilon n,$ $\left\vert
V_{1}\right\vert =...=\left\vert V_{p}\right\vert ;$

\emph{(ii)} for every $i\in\left[  p\right]  ,$ all but at most $\varepsilon
p$ pairs $\left(  V_{i},V_{j}\right)  ,$ $\left(  j\in\left[  p\right]
\backslash\left\{  i\right\}  \right)  ,$ are $\varepsilon$-uniform.
\end{theorem}

Usually SUL is stated with a weaker and less convenient form of condition
\emph{(ii)}; the above form, however, is easily implied.

We present below some counting lemmas needed in the proof of the main theorem.
Lemmas of this kind are known and their proofs are routine, nevertheless, for
the sake of completeness, we present them in some detail.

For every integer $r\geq0,$ set
\[
\phi\left(  r\right)  =r!\sum_{i=0}^{r-1}\frac{1}{i!}.
\]

\begin{lemma}
\label{Xsec}Let $\varepsilon>0,$ $r\geq1,$ and $\left(  A,B\right)  $ be an
$\varepsilon$-uniform pair with $d\left(  A,B\right)  =d.$ If $Y\subset B$ and
$\left(  d-\varepsilon\right)  ^{r-1}\left\vert Y\right\vert >\varepsilon
\left\vert B\right\vert ,$ then fewer than
\[
\varepsilon\phi\left(  r\right)  \binom{\left\vert A\right\vert }{r}%
\]
$r$-sets $R\subset A$ satisfy%
\begin{equation}
d_{Y}\left(  R\right)  \leq\left(  d-\varepsilon\right)  ^{r}\left\vert
Y\right\vert . \label{upc}%
\end{equation}

\end{lemma}

\begin{proof}
Since this result is essentially known (see \cite{KoSi93}, Fact 1.4), we shall
only sketch the proof. We use induction on $r.$ Let $\mathcal{F}_{r}$ be the
class of $r$-sets in $A$ satisfying inequality (\ref{upc}). Observe that those
members of $\mathcal{F}_{r+1}$ that contain a member of $\mathcal{F}_{r}$ are
at most $\left\vert \mathcal{F}_{r}\right\vert \left(  \left\vert A\right\vert
-r\right)  n;$ also, for every $r$-set $R\notin\mathcal{F}_{r},$ at most
$\varepsilon\left\vert A\right\vert $ members of $\mathcal{F}_{r+1}$ contain
$R.$ Therefore,%
\[
\left\vert \mathcal{F}_{r+1}\right\vert \leq\varepsilon\left\vert A\right\vert
\binom{\left\vert A\right\vert -1}{r}+\left\vert \mathcal{F}\right\vert
\left(  \left\vert A\right\vert -r\right)
\]
and the assertion follows.
\end{proof}

With a simple change of signs we obtain a twin result.

\begin{lemma}
\label{Xsec1}Let $\varepsilon>0,$ $r\geq1,$ and $\left(  A,B\right)  $ be an
$\varepsilon$-uniform pair with $d\left(  A,B\right)  =d.$ If $Y\subset B$ and
$\left(  d+\varepsilon\right)  ^{r-1}\left\vert Y\right\vert >\varepsilon
\left\vert B\right\vert ,$ then fewer than
\[
\varepsilon\phi\left(  r\right)  \binom{\left\vert A\right\vert }{r}%
\]
$r$-sets $R\subset A$ satisfy%
\begin{equation}
d_{Y}\left(  R\right)  \geq\left(  d+\varepsilon\right)  ^{r}\left\vert
Y\right\vert . \label{upc1}%
\end{equation}

\end{lemma}

Lemmas \ref{Xsec} and \ref{Xsec1} imply the following statement.

\begin{lemma}
\label{Xsec2}Let $\varepsilon>0,$ $r\geq1,$ and $\left(  A,B\right)  $ be an
$\varepsilon$-uniform pair with $d\left(  A,B\right)  =d.$ Then:

\emph{(i)} at least
\[
\left(  1-\varepsilon\phi\left(  r\right)  \right)  \binom{\left\vert
A\right\vert }{r}%
\]
$r$-sets $R\subset A$ satisfy
\[
\widehat{d}_{B}\left(  R\right)  -d^{r}\left\vert B\right\vert >-\varepsilon
r\left\vert B\right\vert ;
\]

\emph{(ii)} at least
\[
\left(  1-\varepsilon\phi\left(  r\right)  \right)  \binom{\left\vert
A\right\vert }{r}%
\]
$r$-sets $R\subset A$ satisfy%
\[
\widehat{d}_{B}\left(  R\right)  -d^{r}\left\vert B\right\vert <\varepsilon
r\left\vert B\right\vert .
\]

\end{lemma}

\begin{proof}
To prove assertion \emph{(i)}, observe first that it holds trivially if
$d^{r}<\varepsilon r.$ On the other hand, $d^{r}\geq\varepsilon r$ implies
$\left(  d-\varepsilon\right)  ^{r-1}\geq\varepsilon;$ applying Lemma
\ref{Xsec1} with $Y=B,$ we deduce that at least
\[
\left(  1-\varepsilon\phi\left(  r\right)  \right)  \binom{\left\vert
A\right\vert }{r}%
\]
$r$-sets $R\subset A$ satisfy
\[
d_{B}\left(  R\right)  >\left(  d-\varepsilon\right)  ^{r}\left\vert
B\right\vert >\left(  d^{r}-r\varepsilon d^{r-1}\right)  \left\vert
B\right\vert \geq\left(  d^{r}-r\varepsilon\right)  \left\vert B\right\vert ,
\]
completing the proof of \emph{(i)}.

To prove assertion \emph{(ii)}, we use induction on $r.$ For $r=1$ the
assertion follows from Lemma \ref{Xsec1} with $Y=B;$ assume $r\geq2$ and the
assertion true for $r^{\prime}<r.$

Observe that, if $\varepsilon>1-d,$ we deduce
\[
d^{r}+\varepsilon r>d^{r}+\left(  1-d\right)  r\geq1,
\]
and the assertion follows trivially. From $\varepsilon\leq1-d$ we find that
\[
d^{r}+\varepsilon r\geq d^{r}+\varepsilon\left(  \left(  d+\varepsilon\right)
^{r-1}+\varepsilon\left(  d+\varepsilon\right)  ^{r-2}+...+\varepsilon
^{r-1}\right)  =\left(  d+\varepsilon\right)  ^{r}%
\]
so, provided $\left(  d+\varepsilon\right)  ^{r-1}\geq\varepsilon$ holds, we
may apply Lemma \ref{Xsec1} with $Y=B$ and complete the proof of \emph{(ii)}.

It remains to consider the case $\left(  d+\varepsilon\right)  ^{r-1}%
<\varepsilon$ which is only possible if $r>2.$ Let $\mathcal{F}$ be the family
of all $\left(  r-1\right)  $-sets $R\subset A$ satisfying%
\[
\widehat{d}_{B}\left(  R\right)  -d^{r-1}\left\vert B\right\vert
<\varepsilon\left(  r-1\right)  \left\vert B\right\vert ;
\]
by the inductive assumption,
\[
\left\vert \mathcal{F}\right\vert >\left(  1-\varepsilon\phi\left(
r-1\right)  \right)  \binom{\left\vert A\right\vert }{r-1}.
\]
If an $r$-set $R\subset A$ contains a member $R^{\prime}\in\mathcal{F},$ we
find that
\[
\widehat{d}_{B}\left(  R\right)  \leq\widehat{d}_{B}\left(  R^{\prime}\right)
<\left(  \varepsilon\left(  r-1\right)  +d^{r-1}\right)  \left\vert
B\right\vert <\varepsilon r\left\vert B\right\vert \leq\left(  \varepsilon
r+d^{r}\right)  \left\vert B\right\vert .
\]
Since there are at least%
\[
\frac{\left\vert \mathcal{F}\right\vert \left(  n-r+1\right)  }{r}>\left(
1-\varepsilon\phi\left(  r-1\right)  \right)  \binom{\left\vert A\right\vert
}{r}>\left(  1-\varepsilon\phi\left(  r\right)  \right)  \binom{\left\vert
A\right\vert }{r}%
\]
such $r$-sets, the proof is completed.
\end{proof}

Next we shall present a similar result for pairs across different vertex classes.

\begin{lemma}
\label{XPle2} Let $\varepsilon>0$ and $\left(  A_{1},B\right)  ,$ $\left(
A_{2},B\right)  $ be $\varepsilon$-uniform pairs with $d\left(  A_{1}%
,B\right)  =d_{1}$ and\ \ $d\left(  A_{2},B\right)  =d_{2}.$ Then:

\emph{(i) }at least $\left(  1-2\varepsilon\right)  \left\vert A_{1}%
\right\vert \left\vert A_{2}\right\vert $ pairs $\left(  u,v\right)  \in
A_{1}\times A_{2}$ satisfy
\[
\widehat{d}_{B}\left(  uv\right)  -d_{1}d_{2}\left\vert B\right\vert
>-2\varepsilon\left\vert B\right\vert ;
\]

\emph{(ii)} at least $\left(  1-2\varepsilon\right)  \left\vert A_{1}%
\right\vert \left\vert A_{2}\right\vert $ pairs $\left(  u,v\right)  \in
A_{1}\times A_{2}$ satisfy
\[
\widehat{d}_{B}\left(  uv\right)  -d_{1}d_{2}\left\vert B\right\vert
<2\varepsilon\left\vert B\right\vert .
\]

\end{lemma}

\begin{proof}
To prove assertion \emph{(i)}, observe first that it holds trivially if
$d_{1}<2\varepsilon$ or $d_{2}<2\varepsilon,$ so we shall assume $d_{1}%
\geq2\varepsilon$ and\ $d_{2}\geq2\varepsilon.$ Let%
\[
X=\left\{  u\in A_{1}:d_{B}\left(  u\right)  \leq\left(  d_{1}-\varepsilon
\right)  \left\vert B\right\vert \right\}  .
\]
Applying Lemma \ref{Xsec} to the pair $\left(  A_{1},B\right)  $ with $r=1,$
$Y=B$, we find that $\left\vert X\right\vert <\varepsilon\left\vert
A_{1}\right\vert .$ Select any $u\in A_{1}\backslash X,$ and let
\[
Y=\left\{  v\in A_{2}:\widehat{d}_{B}\left(  uv\right)  \leq\left(
d_{2}-\varepsilon\right)  d_{B}\left(  u\right)  \right\}  .
\]
Applying Lemma \ref{Xsec} to the pair $\left(  A_{2},B\right)  $ with $r=1$
and $Y=\Gamma\left(  u\right)  \cap B$ , we find that $\left\vert Y\right\vert
<\varepsilon\left\vert A_{2}\right\vert .$ Therefore, at least%
\[
\left(  1-\varepsilon\right)  ^{2}\left\vert A_{1}\right\vert \left\vert
A_{2}\right\vert >\left(  1-2\varepsilon\right)  \left\vert A_{1}\right\vert
\left\vert A_{2}\right\vert
\]
pairs $\left(  u,v\right)  \in A_{1}\times A_{2}$ satisfy
\[
\widehat{d}_{B}\left(  uv\right)  >\left(  d_{1}-\varepsilon\right)  \left(
d_{2}-\varepsilon\right)  \left\vert B\right\vert >d_{1}d_{2}\left\vert
B\right\vert -2\varepsilon\left\vert B\right\vert ,
\]
completing the proof of \emph{(i)}.

To prove assertion \emph{(ii)}, observe first that, if
\begin{equation}
\left(  d_{1}+\varepsilon\right)  \left(  d_{2}+\varepsilon\right)
>d_{1}d_{2}+2\varepsilon, \label{in3}%
\end{equation}
we deduce
\[
d_{1}d_{2}+2\varepsilon>4-2d_{1}-2d_{2}+d_{1}d_{2}\geq\left(  2-d_{1}\right)
\left(  2-d_{2}\right)  \geq1,
\]
and the assertion follows trivially, so we shall assume that (\ref{in3})
fails. Applying the same argument as in the proof of \emph{(i)}, we find that
at least $\left(  1-2\varepsilon\right)  \left\vert A_{1}\right\vert
\left\vert A_{2}\right\vert $ pairs $\left(  u,v\right)  \in A_{1}\times
A_{2}$ satisfy the inequality
\[
\widehat{d}_{B}\left(  uv\right)  <\left(  d_{1}+\varepsilon\right)  \left(
d_{2}+\varepsilon\right)  \left\vert B\right\vert \leq\left(  d_{1}%
d_{2}+2\varepsilon\right)  \left\vert B\right\vert ,
\]
completing the proof of \emph{(ii)}.
\end{proof}

\subsection{Sums and averages of codegrees}

In this subsection we shall investigate codegrees in graphs consisting of
several $\varepsilon$-uniform pairs.

\begin{lemma}
\label{dle} Let $\varepsilon>0$ and $H$ be a graph whose vertices are
partitioned as
\[
V\left(  H\right)  =A\cup B_{1}\cup...\cup B_{p}%
\]
so that
\[
\left\vert A\right\vert =\left\vert B_{1}\right\vert =...=\left\vert
B_{p}\right\vert =t.
\]
For every $i\in\left[  p\right]  ,$ let the pair $\left(  A,B_{i}\right)  $ be
$\varepsilon$-uniform and set $d\left(  A,B_{i}\right)  =d_{i}$. Then the
inequality%
\[
\left\vert \sum_{\left\{  u,v\right\}  \in S}\sum_{i=1}^{p}\widehat{d}_{B_{i}%
}\left(  uv\right)  -t\left\vert S\right\vert \sum_{i=1}^{p}d_{i}%
^{2}\right\vert <5p\varepsilon t^{3}%
\]
holds for every set $S$ of $2$-sets in $A.$
\end{lemma}

\begin{proof}
We shall prove first that, for every $i\in\left[  p\right]  ,$
\begin{equation}
-5\varepsilon t^{3}\leq\sum_{\left\{  u,v\right\}  \in S}\widehat{d}_{B_{i}%
}\left(  uv\right)  -t\left\vert S\right\vert d_{i}^{2}\leq5\varepsilon t^{3}.
\label{ineq1}%
\end{equation}

Indeed, applying Lemma \ref{Xsec2} to the pair $\left(  A,B_{i}\right)  $ with
$r=2$ and $Y=B_{i},$ we find that at least $\left\vert S\right\vert
-4\varepsilon t^{2}$ sets $\left\{  u,v\right\}  \in S$ satisfy
\[
-2\varepsilon t<\widehat{d}_{B_{i}}\left(  uv\right)  -d_{i}^{2}t<2\varepsilon
t,
\]
and, therefore,%
\[
-2\varepsilon t\left\vert S\right\vert -4\varepsilon t^{3}<\sum_{\left\{
u,v\right\}  \in S}\widehat{d}_{B_{i}}\left(  uv\right)  -d_{i}^{2}t\left\vert
S\right\vert <2\varepsilon t\left\vert S\right\vert +4\varepsilon t^{3}.
\]

As $\left\vert S\right\vert <t^{2}/2$, inequality (\ref{ineq1}) follows;
summing it for $i=1,...,p$ we obtain the desired result.
\end{proof}

\begin{corollary}
\label{cdle} Under the conditions of Lemma \ref{dle}, if $\left\vert
S\right\vert \geq\alpha t^{2}$ for some $\alpha>0,$ then,%
\[
\left\vert \frac{1}{\left\vert S\right\vert }\sum_{\left\{  u,v\right\}  \in
S}\sum_{i=1}^{p}\widehat{d}_{B_{i}}\left(  uv\right)  -t\sum_{i=1}^{p}%
d_{i}^{2}\right\vert <\frac{5p\varepsilon}{\alpha}t.
\]

\end{corollary}

\begin{lemma}
\label{lebs} Suppose $\varepsilon$ $>0$ and $H$ is a graph whose vertices are
partitioned as
\[
V\left(  H\right)  =A_{1}\cup A_{2}\cup B_{1}\cup...\cup B_{p}%
\]
so that
\[
\left\vert A_{1}\right\vert =\left\vert A_{2}\right\vert =\left\vert
B_{1}\right\vert =...=\left\vert B_{k}\right\vert =t.
\]
For every $i\in\left[  2\right]  ,$ $j\in\left[  k\right]  ,$ let the pair
$\left(  A_{i},B_{j}\right)  $ be $\varepsilon$-uniform and set $d\left(
A_{i},B_{j}\right)  =d_{ij}$. Then, the inequality
\[
\left\vert \sum_{\left(  u,v\right)  \in S}\sum_{i=1}^{p}\widehat{d}_{B_{i}%
}\left(  uv\right)  -t\left\vert S\right\vert \sum_{i=1}^{k}d_{1i}%
d_{2i}\right\vert <6\varepsilon pt^{3}%
\]
holds for every set $S\subset A_{1}\times A_{2}.$
\end{lemma}

\begin{proof}
We shall prove first that, for every $i\in\left[  p\right]  ,$
\begin{equation}
-6\varepsilon t^{3}\leq\sum_{\left(  u,v\right)  \in S}\widehat{d}_{B_{i}%
}\left(  uv\right)  -t\left\vert S\right\vert d_{1i}d_{2i}\leq6\varepsilon
t^{3}. \label{ineq2}%
\end{equation}

Indeed, applying Lemma \ref{XPle2} with $B=B_{i},$ we find that at least
$\left\vert S\right\vert -4\varepsilon t^{2}$ pairs $\left(  u,v\right)  \in
S$ satisfy
\[
-2\varepsilon t<\widehat{d}_{B_{i}}\left(  uv\right)  -d_{1i}d_{2i}%
t<2\varepsilon t,
\]
and, therefore,%
\[
-2\varepsilon t\left\vert S\right\vert -4\varepsilon t^{3}<\sum_{\left(
u,v\right)  \in S}\widehat{d}_{B_{i}}\left(  uv\right)  -d_{1i}d_{2i}%
t\left\vert S\right\vert <2\varepsilon t\left\vert S\right\vert +4\varepsilon
t^{3}.
\]

As $\left\vert S\right\vert \leq t^{2}$, inequality (\ref{ineq2}) follows;
summing it for $i=1,...,p$ we obtain the desired result.
\end{proof}

\begin{corollary}
\label{clebs} Under the conditions of Lemma \ref{lebs}, if $\left\vert
S\right\vert \geq\alpha t^{2}$ for some $\alpha>0,$ then,%
\[
\left\vert \frac{1}{\left\vert S\right\vert }\sum_{\left(  u,v\right)  \in
S}\sum_{i=1}^{p}\widehat{d}_{B_{i}}\left(  uv\right)  -t\sum_{i=1}^{p}%
d_{1i}d_{2i}\right\vert <\frac{6p\varepsilon}{\alpha}t.
\]

\end{corollary}

\section{\label{prf}Proof of the main theorem}

Let $d\geq0,$ $a\geq0,$ $c\geq0.$ A sequence $\left\{  SR\left(  n_{s}%
,k_{s},\lambda_{s},\mu_{s}\right)  \right\}  _{s=1}^{\infty}$ of srgs of
increasing order such that
\[
\lim_{s\rightarrow\infty}\frac{k_{s}}{n_{s}}=d,\text{ \ \ }\lim_{s\rightarrow
\infty}\frac{\lambda_{s}}{n_{s}}=a,\text{ \ \ }\lim_{s\rightarrow\infty}%
\frac{\mu_{s}}{n_{s}}=c
\]
is called a $CSR\left(  d,a,c\right)  $\emph{ sequence}.

Note that to prove Theorem \ref{main} it suffices to show that the parameters
$d,a,c$ of any $CSR\left(  d,a,c\right)  $ sequence of nontrivial srgs satisfy
the equalities%
\begin{equation}
a=c=d^{2}. \label{meq}%
\end{equation}

Indeed, assume Theorem \ref{main} false - that is to say, there exist
$\varepsilon>0$ and a sequence $\left\{  SR\left(  n_{s},k_{s},\lambda_{s}%
,\mu_{s}\right)  \right\}  _{s=1}^{\infty}$ of nontrivial srgs of increasing
order such that
\begin{equation}
\left\vert \frac{\lambda_{s}}{n_{s}}-\frac{k_{s}^{2}}{n_{s}^{2}}\right\vert
>\varepsilon\text{ \ \ or \ \ }\left\vert \frac{\mu_{s}}{n_{s}}-\frac
{k_{s}^{2}}{n_{s}^{2}}\right\vert >\varepsilon. \label{ineq3}%
\end{equation}
From the sequence $\left\{  SR\left(  n_{s},k_{s},\lambda_{s},\mu_{s}\right)
\right\}  _{s=1}^{\infty}$ we can always select a subsequence that is a
$CSR\left(  d,a,c\right)  $ sequence for some $d\geq0,$ $a\geq0,$ $c\geq0;$ in
view of inequalities (\ref{ineq3}), condition (\ref{meq}) fails, as claimed.

To prove equalities (\ref{meq}) we shall establish some facts about
$CSR\left(  d,a,c\right)  $ sequences. Observe first that, if $\left\{
G_{s}\right\}  _{s=1}^{\infty}$ is a $CSR\left(  d,a,c\right)  $ sequence,
then $\left\{  \overline{G_{s}}\right\}  _{s=1}^{\infty}$ is a
\[
CSR\left(  1-d,1-2d+c,1-2d+a\right)
\]
sequence.

Also, the well-known relations%
\[
k>\lambda\text{, \ \ }k\geq\mu,\text{ \ \ }k\left(  k-\lambda-1\right)
=\left(  n-k-1\right)  \mu,
\]
holding for any $SR\left(  n,k,\lambda,\mu\right)  ,$ imply that the
parameters of any $CSR\left(  d,a,c\right)  $ sequence satisfy
\begin{equation}
d\geq a,\text{ \ \ }d\geq c, \label{eq0}%
\end{equation}

\begin{equation}
d^{2}-\left(  a-c\right)  d-c=0. \label{eq1}%
\end{equation}

Thus, equalities (\ref{meq}) hold for $d=0,$ and, applying the same argument
to $\left\{  \overline{G_{s}}\right\}  _{s=1}^{\infty},$ they hold for $d=1$
as well. Therefore, we may and shall assume that $0<d<1.$

\begin{lemma}
\label{ltr} If $0<d<1$ and $\left\{  G_{s}\right\}  _{s=1}^{\infty}$ is a
$CSR\left(  d,a,c\right)  $ sequence of nontrivial srgs then $d\neq a$ and
$d\neq c.$
\end{lemma}

\begin{proof}
Assume $d=a;$ then equality (\ref{eq1}) implies $c=0.$ We shall show that
$p=d^{-1}$ is integer and for $s$ sufficiently large, $G_{s}$ is a union of
$p$ complete graphs of equal order.

Let $n_{s},k_{s},\lambda_{s},\mu_{s}$ be the parameters of $G_{s}$ for
$s=1,2,....$ Select any $u\in V\left(  G_{s}\right)  $ and let $\Gamma\left(
u\right)  $ be the set of its neighbors. Clearly, $\left\vert \Gamma\left(
u\right)  \right\vert =k_{s}$ and the graph $G_{s}\left[  \Gamma\left(
u\right)  \right]  $ is $\lambda_{s}$-regular. If $v,w\in\Gamma\left(
u\right)  $ are two nonadjacent vertices, then, by the inclusion-exclusion
formula, we find that
\[
\widehat{d}_{\Gamma\left(  u\right)  }\left(  vw\right)  \geq2\lambda
_{s}-k_{s}=k_{s}+o\left(  n_{s}\right)  ,
\]
and hence $c=d.$ Thus $d=0,$ a contradiction. We conclude that $G\left[
\Gamma\left(  u\right)  \right]  $ is a complete graph of order $k_{s}.$

Furthermore, $\Gamma\left(  u\right)  \cap\Gamma\left(  v\right)
=\varnothing$ for any two nonadjacent vertices $u,v\in V\left(  G_{s}\right)
.$ Indeed, if $w\in\Gamma\left(  u\right)  \cap\Gamma\left(  v\right)  ,$ then
$u,v\in\Gamma\left(  w\right)  ,$ and, therefore, must be adjacent, contrary
to our choice. Thus for any $u\in V\left(  G_{s}\right)  ,$ the set
$\Gamma\left(  u_{i}\right)  \cup\left\{  u_{i}\right\}  $ is a complete graph
of order $k_{s}+1$.

Select a maximal independent set $\left\{  u_{1},...,u_{p}\right\}  $ in
$G_{s}.$ Since $\left\{  u_{1},...,u_{p}\right\}  $ is maximal, we have
\[
\cup_{i=1}^{p}\left(  \Gamma\left(  u_{i}\right)  \cup\left\{  u_{i}\right\}
\right)  =V\left(  G_{s}\right)  .
\]
Thus $d=1/p$ and $V\left(  G_{s}\right)  $ is partitioned in $p$ complete
graphs of order $k_{s}+1$. To complete the proof we have to show that no edge
joins vertices from different complete graphs.

Let $uv$ be an edge such that $u\in\Gamma\left(  u_{i}\right)  \cup\left\{
u_{i}\right\}  ,$ $v\in\Gamma\left(  u_{j}\right)  \cup\left\{  u_{j}\right\}
$, and $i\neq j.$ Since $\Gamma\left(  v\right)  $ is a complete graph and
$u\in\Gamma\left(  v\right)  ,$ then $u$ is adjacent to all vertices of
$\Gamma\left(  u_{j}\right)  \cup\left\{  u_{j}\right\}  ,$ implying $d\left(
u\right)  \geq2k_{s}+1,$ a contradiction, completing the proof.

The case $d=c$ follows by applying the above argument to the sequence
$\left\{  \overline{G_{s}}\right\}  _{s=1}^{\infty}.$
\end{proof}

\begin{proof}
[Proof of Theorem \ref{main}]Let $\left\{  G_{s}\right\}  _{s=1}^{\infty}$ be
a $CSR\left(  d,a,c\right)  $ sequence of nontrivial srgs and suppose
$n_{s},k_{s},\lambda_{s},\mu_{s}$ are the parameters of $G_{s}$ for
$s=1,2,...$. Our goal is to prove equalities (\ref{meq}). Note that, it
suffices to prove that $a=c,$ for, then, the equality $a=d^{2}$ follows
immediately from equality (\ref{eq1}). Observe that since $G_{s}$ are
nontrivial, by Lemma \ref{ltr} we have%
\[
d\neq a,\ \ d\neq c.
\]

Assume
\[
a\neq c,
\]
set
\begin{equation}
\delta=\min\left\{  \left\vert a-c\right\vert ,\left\vert d-a\right\vert
,\left\vert d-c\right\vert ,\frac{1}{10}\right\}  , \label{delt}%
\end{equation}
and let
\begin{align*}
\varepsilon &  =\left(  \frac{\delta}{20}\right)  ^{2},\\
l  &  =\left\lceil 1/\varepsilon\right\rceil .
\end{align*}

Select $s$ so large that the inequalities
\begin{align}
\left\vert k_{s}-dn_{s}\right\vert  &  <\varepsilon n_{s},\label{mins1}\\
\left\vert \lambda_{s}-an_{s}\right\vert  &  <\varepsilon n_{s},\label{mins2}%
\\
\left\vert \mu_{s}-cn_{s}\right\vert  &  <\varepsilon n_{s}\nonumber
\end{align}
hold and, in addition, $n_{s}$ is large enough to apply SUL to $G_{s}$ with
parameters $\varepsilon$ and $l;$ for technical reasons we also require that
$n_{s}>3M\left(  \varepsilon,l\right)  .$

Thus there is a partition $V\left(  G_{s}\right)  =\cup_{i=0}^{p}V_{i}$ such
that $l\leq p\leq M\left(  \varepsilon,l\right)  $ and:

\emph{i)} $\left\vert V_{0}\right\vert <\varepsilon\left\vert G_{s}\right\vert
,$ $\left\vert V_{1}\right\vert =...=\left\vert V_{p}\right\vert ;$

\emph{ii)} for every $i\in\left[  p\right]  ,$ all but at most $\varepsilon p$
pairs $\left(  V_{i},V_{j}\right)  ,$ $\left(  j\in\left[  p\right]
\backslash\left\{  i\right\}  \right)  ,$ are $\varepsilon$-uniform.

Let $n=n_{s},$ $t=\left\vert V_{1}\right\vert ,$ and set $d_{ij}=d\left(
V_{i},V_{j}\right)  $ for every $i,j\in\left[  p\right]  ,$ $\left(  i\neq
j\right)  $. Observe that the inequality $n>3M\left(  \varepsilon,l\right)  $
and condition \emph{(i)} imply%
\begin{equation}
2\leq t\leq\frac{n}{p}\leq\frac{n}{l}\leq\varepsilon n \label{maxt}%
\end{equation}
and
\begin{equation}
tp\leq n\leq\frac{tp}{1-\varepsilon}<tp\left(  1+2\varepsilon\right)  .
\label{maxn}%
\end{equation}

Our first goal is to prove that, if the inequalities
\begin{equation}
\sqrt{\varepsilon}t^{2}<e\left(  V_{i},V_{j}\right)  <\left(  1-\sqrt
{\varepsilon}\right)  t^{2} \label{parc}%
\end{equation}
hold for some pair $\left(  V_{i},V_{j}\right)  ,$ then the inequality%
\begin{equation}
\left\vert a-c\right\vert <\delta, \label{in1}%
\end{equation}
holds, contradicting the choice of $\delta$.

Suppose a pair $\left(  V_{i},V_{j}\right)  $ satisfies inequalities
(\ref{parc}). Let
\[
R=\left\{  r:r\in\left[  p\right]  \backslash\left\{  i,j\right\}  ,\text{
}\left(  V_{i},V_{r}\right)  \text{ and }\left(  V_{j},V_{r}\right)  \text{
are }\varepsilon\text{-uniform}\right\}  .
\]
Observe first that condition \emph{(ii)} implies $\left\vert R\right\vert
\geq\left(  1-2\varepsilon\right)  p.$ Select any vertex $u\in V_{i};$
inequality (\ref{mins1}) implies
\[
\left(  d-\varepsilon\right)  n<d\left(  u\right)  <\left(  d+\varepsilon
\right)  n,
\]
and, therefore,%
\[
\left(  d-\varepsilon\right)  n<\sum_{r=0}^{p}d_{V_{r}}\left(  u\right)
<\left(  d+\varepsilon\right)  n.
\]
Hence, in view of $\left\vert R\right\vert \geq\left(  1-2\varepsilon\right)
p$ and $pt\leq n,$ we deduce%
\[
\left(  d-4\varepsilon\right)  n<\left(  d-\varepsilon\right)  n-2\varepsilon
pt<\sum_{r\in R}d_{V_{r}}\left(  u\right)  <\left(  d+\varepsilon\right)  n,
\]
and, by inequalities (\ref{maxn}), it follows that%
\[
\left(  d-4\varepsilon\right)  pt<\sum_{r\in R}d_{V_{r}}\left(  u\right)
<\left(  d+\varepsilon\right)  \left(  1+2\varepsilon\right)  pt\leq\left(
d+4\varepsilon\right)  pt.
\]
Summing this inequality for all $u\in V_{i}$ and dividing by $t^{2},$ we
obtain
\begin{equation}
\left(  d-4\varepsilon\right)  p<\sum_{r\in R}d_{ir}<\left(  d+4\varepsilon
\right)  p; \label{deq1}%
\end{equation}
by symmetry we also have%
\begin{equation}
\left(  d-4\varepsilon\right)  p<\sum_{r\in R}d_{jr}<\left(  d+4\varepsilon
\right)  p. \label{deq2}%
\end{equation}

Applying Corollary \ref{clebs} with $A_{1}=V_{i},$ $A_{2}=V_{j},$ $B_{r}%
=V_{r}$ for all $r\in R,$ and $S=E\left(  V_{1},V_{2}\right)  ,$ we see that
\begin{equation}
\left\vert \frac{1}{e\left(  V_{i},V_{j}\right)  }\sum_{\left(  u,v\right)
\in E\left(  V_{1},V_{2}\right)  }\frac{1}{t}\sum_{r\in R}\widehat{d}_{V_{r}%
}\left(  uv\right)  -\sum_{r\in R}d_{1r}d_{2r}\right\vert <\frac
{6p\varepsilon}{\sqrt{\varepsilon}}=6\sqrt{\varepsilon}p. \label{in2}%
\end{equation}

Furthermore, select any edge $uv$ such that $u\in V_{i}$ and $v\in V_{j}.$
Condition (\ref{mins2}) implies
\[
\left(  a-\varepsilon\right)  n<\widehat{d}\left(  uv\right)  <\left(
a+\varepsilon\right)  n;
\]
conditions \emph{(i)} and \emph{(ii)} imply
\[
0\leq\widehat{d}\left(  uv\right)  -\sum_{r\in R}\widehat{d}_{V_{r}}\left(
uv\right)  =\widehat{d}_{V_{0}}\left(  uv\right)  +\sum_{r\in\left[  p\right]
\backslash R}\widehat{d}_{V_{r}}\left(  uv\right)  <\varepsilon n+2\varepsilon
pt<3\varepsilon n.
\]
It follows that
\[
\left(  a-4\varepsilon\right)  n<\sum_{r\in R}\widehat{d}_{V_{r}}\left(
uv\right)  <\left(  a+\varepsilon\right)  n
\]
and, estimating $n$ from (\ref{maxn}), we see that
\[
\left(  a-4\varepsilon\right)  p<\frac{1}{t}\sum_{r\in R}\widehat{d}_{V_{r}%
}\left(  uv\right)  <\left(  a+\varepsilon\right)  \left(  1+2\varepsilon
\right)  p<\left(  a+4\varepsilon\right)  p.
\]
Hence, inequality (\ref{in2}) implies
\begin{equation}
\left(  a-10\sqrt{\varepsilon}\right)  p<\sum_{r\in R}d_{1r}d_{2r}<\left(
a+10\sqrt{\varepsilon}\right)  p. \label{aeq}%
\end{equation}

Applying the same argument to any pair $\left(  u,v\right)  \in V_{i}\times
V_{j}$ such that $uv\notin E\left(  V_{i},V_{j}\right)  $, we obtain
\[
\left(  c-10\sqrt{\varepsilon}\right)  p<\sum_{r\in R}d_{1r}d_{2r}<\left(
c+10\sqrt{\varepsilon}\right)  p.
\]
These inequalities together with inequalities (\ref{aeq}) imply%
\[
\left\vert a-c\right\vert <20\sqrt{\varepsilon}\leq\delta,
\]
as claimed.

Therefore, we may and shall assume that condition (\ref{parc}) fails for all
pairs $\left(  V_{i},V_{j}\right)  $ - that is to say, for every
$i,j\in\left[  p\right]  ,$ $\left(  i\neq j\right)  ,$ either
\[
d_{ij}\leq\sqrt{\varepsilon}\text{ \ \ or \ \ }d_{ij}\geq1-\sqrt{\varepsilon
}.
\]
A simple calculation shows that then
\begin{equation}
0\leq d_{ij}-d_{ij}^{2}\leq\sqrt{\varepsilon} \label{dens}%
\end{equation}
holds for every $i,j\in\left[  p\right]  ,$ $\left(  i\neq j\right)  .$ We
shall prove that these inequalities imply either
\[
\left\vert d-a\right\vert <\delta\text{ \ \ or \ \ }\left\vert d-c\right\vert
<\delta,
\]
contradicting (\ref{delt}).

Assume $e\left(  V_{1}\right)  \geq t^{2}/5$ and let
\[
R=\left\{  r:r\in\left[  p\right]  \backslash\left\{  i\right\}  ,\text{ the
pair }\left(  V_{1},V_{r}\right)  \text{ is }\varepsilon\text{-uniform}%
\right\}  .
\]
As above we establish
\[
\left(  d-4\varepsilon\right)  p<\sum_{r\in R}d_{1r}<\left(  d+4\varepsilon
\right)  p;
\]
Hence, in view of (\ref{dens}), we obtain%
\begin{equation}
\left(  d-5\sqrt{\varepsilon}\right)  p<\sum_{r\in R}d_{1r}^{2}<\left(
d+5\sqrt{\varepsilon}\right)  p. \label{sqd}%
\end{equation}
Applying Corollary \ref{cdle} with $A=V_{1},$ $B_{r}=V_{r}$ for all $r\in R,$
and $S=E\left(  V_{1}\right)  ,$ we see that
\[
\left\vert \frac{1}{e\left(  V_{1}\right)  }\sum_{uv\in E\left(  V_{1}\right)
}\frac{1}{t}\sum_{r\in R}\widehat{d}_{V_{r}}\left(  uv\right)  -\sum_{r\in
R}d_{1r}^{2}\right\vert <\frac{5p\varepsilon}{1/5}=25\varepsilon p.
\]

For any edge $uv$ induced by $V_{1},$ as above, we establish that
\[
\left(  a-4\varepsilon\right)  pt<\sum_{r\in R}\widehat{d}_{V_{r}}\left(
uv\right)  <\left(  a+\varepsilon\right)  \left(  1+2\varepsilon\right)
pt<\left(  a+4\varepsilon\right)  pt.
\]
Hence, inequality (\ref{sqd}) implies%
\[
\left\vert d-a\right\vert <29\varepsilon+5\sqrt{\varepsilon}<\delta,
\]
as claimed.

Assuming $e\left(  V_{1}\right)  <t^{2}/5,$ from $t\geq2,$ we see that the
graph $\overline{G\left[  V_{1}\right]  }$ induces at least $t^{2}/5$ edges.
Applying Corollary \ref{cdle} with $A=V_{1},$ $B_{r}=V_{r}$ for all $r\in R,$
and $S=E\left(  \overline{G\left[  V_{1}\right]  }\right)  ,$ by the above
argument applied to the members of $S$, we see that
\[
\left\vert d-c\right\vert <\delta,
\]
as claimed. The proof is completed.
\end{proof}

\section{Concluding remark}

Curiously enough, in the proof of Theorem \ref{main} we did not make much use
of the essential feature of SUL - the independence of $M\left(  \varepsilon
,l\right)  $ on $n.$ This fact suggests that a more involved approach exists,
possibly leading to effective bounds on the values%
\[
\left\vert \lambda-k^{2}/n\right\vert \text{ \ \ and \ \ }\left\vert \mu
-k^{2}/n\right\vert .
\]

\textbf{Acknowledgement} The author is grateful to Cecil Rousseau for his kind
and helpful assistance.

\end{document}